\documentstyle[12pt]{article}
\title{ Integrals of Braided  Hopf Algebras  }
\author{ Xijing Guo,   Shouchuan Zhang\\
 Department  of Mathematics,
Hunan University\\ Changsha  410082, \
 P.R. China  }
\date{}

\begin{document}
\newtheorem{Theorem}{\quad Theorem}[section]
\newtheorem{Proposition}[Theorem]{\quad Proposition}
\newtheorem{Definition}[Theorem]{\quad Definition}
\newtheorem{Corollary}[Theorem]{\quad Corollary}
\newtheorem{Lemma}[Theorem]{\quad Lemma}
\newtheorem{Example}[Theorem]{\quad Example}
\maketitle \addtocounter{section}{-1}

\begin {abstract}

The faithful quasi-dual $H^d$ and strict quasi-dual $H^{d'}$ of an
infinite braided Hopf algebra $H$ are introduced and  it is proved
that every strict quasi-dual $H^{d'}$ is  an $H$-Hopf module. The
connection between the integrals and the maximal rational
$H^{d}$-submodule $H^{d rat }$ of $H^{d}$ is found. That is, $H^{d
rat  }\cong \int ^l_{H^d} \otimes H$ is proved. The existence and
uniqueness of integrals for  braided Hopf algebras in the
Yetter-Drinfeld category $(^B_B{\cal YD},C )$ are given.
\end {abstract}

\section {Introduction}

Integrals of Hopf algebras were introduced by Larson and Sweedler in
\cite {LS69}. Their  connection with the maximal rational
$H^*$-module $H^{*rat }$ of $H^*$ was given by Sweedler in \cite
{Sw69b}, i.e.
\begin {eqnarray}\label {e1}
H^{* rat }\cong \int ^l_{H^*} \otimes H \ \ \ \ \ \ \ \hbox {as }
H\hbox {-Hopf modules}
\end {eqnarray}
\noindent Uniqueness of integrals was proved by Sullivan in \cite
{Su71}. The existence of non-zero integrals was given in \cite
[Theorem 5.3.2] {DNR01}. The integrals have proved to be essential
instruments in constructing invariants of surgically presented
3-manifolds or 3-dimensional topological quantum field theories
\cite {Ke99} \cite {Ku96} \cite {Tu94}.

Recently braided tensor categories were introduced by Joyal and
Street \cite {JS86}. Algebraic structures within them, especially,
braided Hopf algebras or ``braided groups" as well  as cross
products and diagrammatic techniques for such algebraic
constructions were studied by Majid in \cite {Ma93} \cite {Ma94}.
See \cite {Ma95a} \cite {Ma95b} for introductions. Many braided
groups are known, including ones obtained by transmutation \cite
{Ma93} from the (co)quasitriangular Hopf
 algebras, and the universal enveloping algebra
of a Lie color algebra, the Nichols algebras \cite {AS02} and the
Lusztig's quantum algebras \cite {Lu93}. Therefore, it is
interesting to extend the  Hopf algebra constructions to the braided
cases. For finite braided Hopf algebras (braided groups) $H$,\ i.e.
braided Hopf algebra  $H$ with a left dual in braided tensor
categories, Bespalov, Kerler and Lyubashenko \cite {BKL00}, and
Takeuchi \cite {Ta99} introduced an integral and proved that the
integral is an invertible object. Moreover, Takeuchi   proved that
the antipode is an isomorphism and  formula (\ref {e1}) holds.

In this paper we study integrals of infinte braided Hopf algebras in
braided tensor categories. A braided Hopf algebra is called an
infinite braided Hopf algebra if it has no left duals (See \cite
{Ta99}). An important example of infinite braided Hopf algebras is
the universal  enveloping algebra of a Lie superalgebra. So the
integrals of infinite braided Hopf algebras should have important
applications in both mathematics and mathematical physics. We
introduce faithful quasi-dual $H^d$ and strict  quasi-dual $H^d$ of
a braided Hopf algebra $H$. We  prove that every  strict quasi-dual
$H^{d'}$ is an $H$-Hopf module. By imitating Larson and Sweedler's
Hopf module construction, we obtain the connection between the
integrals and the maximal rational $H^{d}$-submodule $H^{d rat }$ of
$H^{d}$. That is, we prove $H^{d rat  }\cong \int ^l_{H^d} \otimes
H$ . We give the existence and uniqueness of integrals for some
infinite braided Hopf algebras living in the Yetter-Drinfeld
category $(^B_B{\cal YD},C )$.

This paper was organized  as follows. In section 1, since it is
possible that $Hom (H, I)$ is not an object in ${\cal C}$ for
braidrd Hopf algebra $H$, we introduce  strict (or faithful)
quasi-dual $H^{d'}$,  and prove that every  strict quasi-dual
$H^{d'}$ is an $H$-Hopf module. In section 2, we concentrate on
braided tensor categories consisting of some braided vector spaces.
We prove $H^{d rat  }\cong \int ^l_{H^d} \otimes H$ for an infinite
braided Hopf algebra $H$ and the maximal rational $H^{d}$-submodule
$H^{d rat }$ of $H^{d}$. That is, we obtain the connection between
integrals and  the maximal rational $H^{d}$-module $H^{d rat }$ of
$H^{d}$. In section 3 we give the existence and uniqueness of
integrals for
 infinite braided Hopf algebras living in the Yetter-Drinfeld category
 $(^B_B{\cal YD},C )$.
In section 4 we show the Maschke's theorem for infinite braided Hopf
algebras.

\section {Strict quasi-duals and  Hopf modules of braided Hopf algebras }
In this section we introduce a  faithful quasi-dual $H^d$ and strict
quasi-dual $H^{d'}$ of braided Hopf algebra $H$ and show that
$H^{d'}$ is an $H$-Hopf module. Using the fundamental theorem of
Hopf modules,
 we show the formula similar to (\ref {e1})
\begin {eqnarray*}\label {e3}
H^{d '}\cong (H ^{d '})^{coH} \otimes H \ \ \ \ \ \ \ \hbox {as }
H\hbox {- Hopf modules in } {\cal C}.
\end {eqnarray*}

We first recall some notations. Let  $({\cal C}, \otimes, I,  C )$
be   a braided tensor category, where $I$ is the identity object and
$C$ is the braiding. We also write $W \otimes f $ for $id _W \otimes
f$ and $f \otimes W$ for  $f \otimes id _W.$ Since every braided
tensor category is always equivalent to a strict braided tensor
category by \cite [Theorem 0.1] {ZC01},
 we may view every braided tensor category as
 a strict braided tensor category.

\begin {Definition} \label {1.1}
Let $H$ be a braided Hopf algebra  in ${\cal C}$. If there is an
algebra $N$ in $ {\cal C}$ and a morphism $<,>$  from $N\otimes H$
to $I$ such that $$<,>(m\otimes H)=(<,>\otimes <,>)(N\otimes
C_{N,H}\otimes H)(N\otimes N\otimes \Delta),\ \ \ \
<\eta,H>=\epsilon$$ then $N$ is called a left quasi-dual of $H$.
Moreover, if  for any objects $U, V$ and four morphisms $f: U
\rightarrow V\otimes N$, $f': U \rightarrow V\otimes N$,
 $g: U \rightarrow H\otimes V$, $g': U \rightarrow H\otimes V$ in ${\cal C}$,
$(V \otimes <,>) (f \otimes H) = (V \otimes <,>) (f' \otimes H)$
implies $f = f'$ and $(<, > \otimes V) (N \otimes g)=(<, > \otimes
V) (N \otimes g')$ implies $g=g'$, then $N$ is called a faithful
quasi-dual of $H$ under $<,>$, written as $H^{d}.$  In addition, if
there are  a left ideal, written as   $H^{d'}$,  of $H^{d}$ \   and
two  morphisms: \ \ $\rightharpoonup: \ \ H \otimes H^ {d}
\rightarrow H^{d} $ and $\rho : \ H^{d'} \rightarrow H^{d'} \otimes
H$ \ in {\cal C } such that
$$<,>(\rightharpoonup\otimes H)=<,>(H^d\otimes m)(H^d\otimes C_{H,H})(C_{H,H^d}\otimes
H)$$
$$ (H^d\otimes <,>)(C_{H^d,H^{d'}}\otimes H)(H^d\otimes
\rho)=m$$ and the constraint $\rightharpoonup $ on $ H \otimes
H^{d'}$ is a morphism to $H^{d'}$, \noindent  then $H^{d'}$ is
called  a  strict quasi-dual of $H$.
\end {Definition}
Let $\leftharpoondown = \rightharpoonup (S \otimes H^{d'} )
C_{H^{d}, H}.$ In fact, if $H$ has a left dual $H^*$ in ${\cal C}$,
then $H^*$ is a strict quasi-dual and faithful quasi-dual of $H$
under evaluation $<,>.$

\begin {Lemma} \label {1.2} Let $H^{d}$ be a faithful quasi-dual of $H$ under $<,>$ and
$C_{H,H} = C_{H,H}^{-1}$.
 Then
$C_{U,V}= (C_{V,U}) ^{-1}$,  for $ U, V = H $ or $H^d.$

\end {Lemma}
{\bf Proof.} See that $(H\otimes <,>)(C_{H^d,H}\otimes
H)=(<,>\otimes H)(H^d\otimes C^{-1}_{H,H})=(<,>\otimes H)(H^d\otimes
C_{H,H})=(H\otimes <,>)(C^{-1}_{H^d,H}\otimes H)$. Thus
$C_{H^d,H}=C^{-1}_{H^d,H}$. See that
$C_{H,H^d}=C^{-1}_{H^d,H}C_{H^d,H}C_{H,H^d}=C^{-1}_{H^d,H}C^{-1}_{H^d,H}C_{H,H^d}=C^{-1}_{H^d,H}$
and $(H^d\otimes <,>)(C_{H^d,H^d}\otimes H)=(<,>\otimes
H^d)(H^d\otimes C^{-1}_{H^d,H})=(<,>\otimes H^d)(H^d\otimes
C_{H^d,H})=(H^d\otimes <,>)(C^{-1}_{H^d,H^d}\otimes H)$. Thus
$C_{H^d,H^d}=C^{-1}_{H^d,H^d}$.$\Box$

If $C_{H,H}= C_{H,H}^{-1}$, then  we say that the braiding is
symmetric on $H$. Throughout this section we always assume that the
braiding is symmetric on $H$. For convenience, for $U, V = H $ or
$H^d$ we  denote the braiding $C_{U,V}$ by $C$.
\begin {Lemma} \label {1.3}
$m(H^d\otimes \leftharpoondown)=\leftharpoondown(m\otimes
H)(\rightharpoonup\otimes H^{d'}\otimes H)(C\otimes H^{d'}\otimes
H)(H^d\otimes C\otimes H)(H^d\otimes H^{d'}\otimes C)(H^d\otimes
H^{d'}\otimes\Delta)$
\end {Lemma}
{\bf Proof.} $ <,>(m\otimes H)(H^d\otimes \leftharpoondown \otimes
H)=(<,>\otimes <,>)(H^d\otimes C\otimes H)(H^d\otimes
\leftharpoondown\otimes\Delta)=<,>(<,>\otimes H^{d'}\otimes
m)(H^d\otimes C\otimes C)(H^d\otimes H^{d'}\otimes C\otimes
H)(H^d\otimes H^{d'}\otimes S\otimes \Delta)$ and
$<,>(\leftharpoondown\otimes H)(m\otimes H\otimes
H)(\rightharpoonup\otimes H^{d'}\otimes H\otimes H)(C\otimes
H^{d'}\otimes H\otimes H)(H^d\otimes C\otimes H\otimes H)(H^d\otimes
H^{d'}\otimes C\otimes H)(H^d\otimes H^{d'}\otimes \Delta\otimes H)
=<,>(m\otimes H)(\rightharpoonup\otimes H^{d'}\otimes C)(H^d\otimes
C\otimes S\otimes H)(H^d\otimes H^{d'}\otimes C\otimes H)(H^d\otimes
H^{d'}\otimes \Delta\otimes H) =<,>(C\otimes
H)(\rightharpoonup\otimes H^{d'}\otimes m)(C\otimes H^{d'}\otimes
C)(H^d\otimes C\otimes S\otimes H)(H^d\otimes H^{d'}\otimes C\otimes
H)(H^d\otimes H^{d'}\otimes \Delta\otimes H) =<,>(H^d\otimes
m)(H^d\otimes C)(H^d\otimes H\otimes H\otimes <,>)(H^d\otimes
H\otimes C\otimes H)(H^d\otimes H\otimes H^{d'}\otimes
\Delta)(H^d\otimes H\otimes H^d\otimes m)(H^d\otimes C\otimes
C)(H^d\otimes H^{d'}\otimes C\otimes H)(H^d\otimes H^{d'}\otimes
S\otimes H\otimes H)(H^d\otimes H^{d'}\otimes \Delta\otimes H)
=(<,>\otimes <,>)(H^d\otimes C\otimes H)(H^d\otimes H^{d'}\otimes
m\otimes H)(H^d\otimes H^{d'}\otimes C\otimes H)(H^d\otimes
H^{d'}\otimes H\otimes m\otimes m)(H^d\otimes H^{d'}\otimes H\otimes
H\otimes C\otimes H)(H^d\otimes H^{d'}\otimes H\otimes \Delta\otimes
\Delta)(H^d\otimes H^{d'}\otimes H\otimes C)(H^d\otimes
H^{d'}\otimes C\otimes H)(H^d\otimes H^{d'}\otimes S\otimes H\otimes
H)(H^d\otimes H^{d'}\otimes \Delta\otimes H) =<,>(H^{d'}\otimes
m)(H^d\otimes C\otimes <,>)(H^d\otimes H\otimes C\otimes
H)(H^d\otimes H\otimes H^{d'}\otimes m\otimes m)(H^d\otimes H\otimes
H^{d'}\otimes H\otimes C\otimes H)(H^d\otimes H\otimes H^{d'}\otimes
\Delta\otimes S\otimes S)(H^d\otimes H\otimes H^{d'}\otimes H\otimes
C)(H^d\otimes H\otimes H^{d'}\otimes H\otimes \Delta)(H^d\otimes
C\otimes C)(H^d\otimes H^{d'}\otimes C\otimes H)(H^d\otimes
H^{d'}\otimes \Delta \otimes H) =(<,>\otimes <,>)(H^d\otimes
C\otimes H)(H^d\otimes H^{d'}\otimes m\otimes m)(H^d\otimes
H^{d'}\otimes m\otimes H\otimes C)(H^d\otimes H^{d'}\otimes C\otimes
C\otimes H)(H^d\otimes H^{d'}\otimes H\otimes C\otimes H\otimes
H)(H^d\otimes H^{d'}\otimes S\otimes S\otimes H\otimes H\otimes
H)(H^d\otimes H^{d'}\otimes C\otimes C\otimes H)(H^d\otimes
H^{d'}\otimes \Delta\otimes H\otimes \Delta)(H^d\otimes
H^{d'}\otimes \Delta\otimes H) =(<,>\otimes <,>)(H^d\otimes C\otimes
H)(H^d\otimes H^{d'}\otimes m\otimes m)(H^d\otimes H^{d'}\otimes
H\otimes m\otimes C)(H^d\otimes H^{d'}\otimes C\otimes C\otimes
H)(H^d\otimes H^{d'}\otimes H\otimes C\otimes H\otimes H)(H^d\otimes
H^{d'}\otimes C\otimes H\otimes H\otimes H)(H^d\otimes H^{d'}\otimes
H\otimes S\otimes C\otimes H)(H^d\otimes H^{d'}\otimes S\otimes
\Delta\otimes \Delta)(H^d\otimes H^{d'}\otimes \Delta\otimes H)
=<,>(<,>\otimes H^{d'}m)(H^{d'}\otimes C\otimes C)(H^d\otimes
H^{d'}\otimes C\otimes H)(H^d\otimes H^{d'}\otimes S\otimes
\Delta)$. Thus we complete the proof. $\Box$

\begin {Theorem} \label {1.4}
    $H^{d'}$ is an $H$-Hopf module .
\end {Theorem}
{\bf Proof.} (1)$(H^{d'},\leftharpoondown)$ is a right $H$-module.\\
$<,>(\leftharpoondown\otimes H)(\leftharpoondown\otimes H\otimes H)
=<,>(H^{d'}\otimes m)(H^{d'}\otimes C)(\leftharpoondown\otimes
S\otimes H)=<,>(H^{d'}\otimes m)(H^{d'}\otimes C)(H^{d'}\otimes
H\otimes m)(H^{d'}\otimes H\otimes C)(H^{d'}\otimes S\otimes
S\otimes H)$ and $<,>(\leftharpoondown\otimes H)(H^{d'}\otimes
m\otimes H)=<,>(H^{d'}\otimes C)(H^{d'}\otimes S\otimes
H)(H^{d'}\otimes m\otimes H)=<,>(H^{d'}\otimes m)(H^{d'}\otimes
C)(H^{d'}\otimes m\otimes H)(H^{d'}\otimes C\otimes H)(H^{d'}\otimes
S\otimes S\otimes H)$. Thus
$\leftharpoondown(\leftharpoondown\otimes
H)=\leftharpoondown(H^{d'}\otimes m)$. Obviously,
$\leftharpoondown(H^{d'}\otimes \eta)=id_{H^{d'}}$. Therefore,
$(H^{d'},\leftharpoondown)$ is a right $H$-module.\\
(2) $(H^{d'},\rho)$ is a right $H$-comodule.\\
See that $(H^{d'}\otimes <,>\otimes <,>)(C\otimes C\otimes
H)(H^d\otimes C\otimes H\otimes H)(H^d\otimes H^d\otimes \rho\otimes
H)(H^d\otimes H^d\otimes \rho)=m(H^d\otimes m)=m(m\otimes
H^d)=(H^{d'}\otimes <,>\otimes <,>)(C\otimes C\otimes H)(H^d\otimes
C\otimes \Delta)(H^d\otimes H^d\otimes \rho)$. Thus $(\rho\otimes
H)\rho=(H^{d'}\otimes \Delta)\rho$. We also have that $(id\otimes
\epsilon)\rho=(H^{d'}\otimes <,>)(C\otimes H)(\eta\otimes \rho)
=m(\eta\otimes H^{d'})=id$. Therefore $(H^{d'},\rho)$ is a right $H$-comodule.\\
(3) See that $(H^{d'}\otimes <,>)(C\otimes H)(H^d\otimes
\rho)(H^d\otimes \leftharpoondown)= m(H^d\otimes
\leftharpoondown)=\leftharpoondown(m\otimes
H)(\rightharpoonup\otimes H^{d'}\otimes H)(C\otimes H^{d'}\otimes
H)(H^d\otimes C\otimes H)(H^d\otimes H^{d'}\otimes C)(H^d\otimes
H^{d'}\otimes \Delta)= (<,>\otimes H^{d'})(\rightharpoonup\otimes
H\otimes \leftharpoondown)(C\otimes C\otimes H)(H^{d'}\otimes
C\otimes H\otimes H)(H^d\otimes H^{d'}\otimes C\otimes H)(H^d\otimes
H^{d'}\otimes H\otimes C)(H^d\otimes \rho\otimes \Delta)=
(<,>\otimes H^{d'})(H^d\otimes m\otimes \leftharpoondown)(H^d\otimes
H\otimes C\otimes H)(H^d\otimes C\otimes C)(H^d\otimes \rho\otimes
\Delta)=(<,>\otimes H^{d'})(H^d\otimes C)(H^d\otimes
\leftharpoondown\otimes m)(H^d\otimes H^{d'}\otimes C\otimes
H)(H^d\otimes \rho\otimes \Delta)$. Thus
$\rho\circ\leftharpoondown=(\leftharpoondown\otimes m)(H^{d'}\otimes
C\otimes H)(\rho\otimes \Delta)$. From (1)(2)(3), we complete the
proof. $\Box$

If ${\cal C}$ has equalizers, then the coinvariant $(H^{d'})^{coH}$
of $H$  in $H^{d'}$ is an object in ${\cal C}.$ Here
$(H^{d'})^{coH}$ denotes the equalizer of the diagram

$$H^{d'}   \left.    \begin{array}{c} \rho \\
\longrightarrow \\
\longrightarrow \\
id \otimes \eta
\end{array} \right.  H^{d'} \otimes H \ \ \ .
$$

Combining Theorem \ref {1.4} and the braided Hopf module fundamental
theorem \cite [Theorem 3.4]{Ta99}, we have

\begin {Theorem} \label {1.5} If ${\cal C}$ has
equalizers or $(H^{d'})^{coH}$ is an object in ${\cal C}$, then
   $$ H^{d'} \cong (H^{d'})^{coH}\otimes H \ \ \ \ \  \ \
(\hbox {as } H \hbox {-Hopf modules in } {\cal C}).$$

\end {Theorem}

\section { Connection between integrals and the maximal rational
$H^d$-submodule $H^{d rat}$ of $H^d$}

In this section, we concentrate on braided tensor categories
consisting of some braided vector spaces. We obtain $H^{d rat }\cong
\int ^l_{H^d} \otimes H$ for an infinite  braided Hopf algebra
 $H$ and the maximal rational $H^{d}$-submodule $H^{d rat }$ of $H^{d}$.

Throughout  this section  we  assume the following unless otherwise
stated:  $H$ is a braided Hopf algebra in ${\cal C}$ with
 $C_{H,H}= C^{-1}_{H,H}$ and  $<,>$ is the evaluation of $H$; there is a faithful quasi-dual
$H^{d}\subseteq H^*$
 We also assume that $k$ is a field and there exists a forgetful functor $F: \
\ {\cal C} \rightarrow {}_k{\cal M}$, which is the category of
vector spaces over $k$ such that $F(U \otimes V) = F(U)\otimes F(V)$
and $F(I)= k$.

Now we give the concept of rational $H^{d}$-modules. For
$H^{d}$-module $(M, \alpha )$, if there is a morphism $\rho $ from
$M$ to $M\otimes H$ in ${\cal C}$ such that the condition of
module-comodule compatibility
$$  (MCOM): \ \ (H \otimes <,>)(C \otimes M )(H^d  \otimes \rho ) = \alpha $$
\noindent holds, then $(M, \alpha )$ is called  rational
$H^{d}$-module.

$$M^H = \{ x\in M \mid   h\cdot x = \epsilon (h) x
\hbox { for every } h\in H \}$$ is called the invariant of $H$  on
$M$.  In particular, if $M$ is a regular $H$-module (i.e. the module
operation is $m$ ), then $M^H$  is  written as $\int _H^l.$  We also
denote $$\{ f \in H^* \mid  g * f  =  g(1)f \hbox { for every  }
g\in H^* \}$$ by $\int _{H^*}^l$. Moreover, for some subset $N $ of
$H^*$, we also  denote  $$\{ f \in N \mid  g * f  = g(1)f \hbox {
for any } g\in N \}$$ by $\int _{N}^l$. Every  element in $\int
_{H^*}^l$ is called an integral on $H.$

Dually, if $(M, \phi )$ is a left $H$-comodule, then the set
$$M^{coH} = \{ x \in M \mid \phi (x) =1  \otimes   x
\}$$ is called the coinvariant of $H$  in $M$.

\begin {Corollary} \label {2.2}
Assume ${\cal C}$ has equalizers and there exists   the maximal
rational $H^{d}$-submodule $H^{d rat }$ of regular module $ H^{d}$.
If \  $\rightharpoonup$ is a morphism from $ H \otimes H^{d }$ to $
H^{d}$  and the constraint on   $H \otimes H^{d rat }$ is a morphism
to $H^{d rat }$ in ${\cal C},$ then
   $$ H^{d rat} \cong \int _{H^d}^l\otimes H \ \ \ \
\  \ \ (\hbox {as } H \hbox {-Hopf modules in } {\cal C})\ .$$
\end {Corollary}

{\bf Proof.} For convenience, let $H^{\Box}$ denote $H^{d rat}$.
Obviously, $H^{\Box}$ is a strict quasi-dual of $H$. By Theorem \ref
{1.5}. It suffices to show $\int_{H^{\Box}}^l = (H^{\Box})^{coH}$.
Obviously, $\int _{H^d}^l \subseteq (H^\Box) ^{co H}.$

Conversely, we see $m=((H^\Box) ^{co H}\otimes <,>)(C\otimes
H)(H^d\otimes \rho)=((H^\Box) ^{co H}\otimes <,>)(C\otimes
\eta)=\epsilon \otimes id_{(H^\Box) ^{co H}}$. Thus $(H^\Box) ^{co
H}\subseteq\int _{H^d} ^l$

Consequently, $\int _{H^d} ^l = (H^{\Box})^{CoH}.$ $\Box$

The above corollary is a generalization of Sweedler's relation (\ref
{e1}). In fact, we have
\begin {Corollary} \label {2.3}  If $H $ is an ordinary Hopf algebra,  then
$\int _{H^{* rat }}^l = \int _{H^*}^l$.
\end {Corollary}
{Proof.} Obviously $H^*$ is a faithful quasi-dual of $H$ and $H^{*
rat}$ is a strict quasi-dual of $H$. By Corollary \ref {2.2}, we can
complete the proof. $\Box$

\section {Existence and uniqueness of integrals for Yetter-Drinfeld module
categories }

In this section we  give the existence and uniqueness  of integrals
for braided Hopf algebras in the Yetter-Drinfeld module category
$(^B_B{\cal YD}, C)$. Throughout this section, $H$ is a braided Hopf
algebra in $(^B_B{\cal YD}, C)$ with finite-dimensional Hopf algebra
$B$. Let $b_B$ denote the coevaluation of $B$ and  $\tau : U \otimes
V \rightarrow V\otimes U $ denote the flip
 $\tau (x\otimes y) = y \otimes x.$ If $(M, \alpha )$ is a left $B$-module, we can define a left $B$-module structure
$\alpha _{M^*}$ on $M^* = Hom _k (M, k)$ such that $(b \cdot x^*)(x)
= x^*(S(b)\cdot x)$ for any $b \in B, x \in M, x^* \in M^*$. If $(M,
\phi )$ is a left $B$-comodule,  we can also  define a left
$B$-comodule structure $\phi _{M^*}$ on $M^*$ such that $  (B
\otimes <,>) (\phi _{M^*} \otimes M) =  (S^{-1} \otimes <,>) (\tau
\otimes M) (M^* \otimes \phi )$. In fact, $\phi _{M^*} = (S^{-1}
\otimes \hat \alpha ) (b _B \otimes M^*)$, where  $ <,> ( \hat
\alpha  \otimes M) = (<,> \otimes <,>) (B^* \otimes \tau \otimes M)
(B^* \otimes M^* \otimes \phi ). $

\begin {Lemma} \label {3.1}
(i) If $(M, \alpha , \phi ) \in (^B_B{\cal YD}, C)$, then
 $(M^*, \alpha _{M^*}, \phi _{M^*}) \in (^B_B{\cal YD}, C)$ and the evaluation $<,>$ is a morphism in $(^B_B{\cal YD}, C)$

(ii) If  $(H, \alpha , \phi )$ is a braided Hopf algebra in  $
(^B_B{\cal YD}, C)$ and the antipode $S$ of $B$  satisfies $ S=
S^{-1}$, then $\rightharpoonup $ is a morphism from $H \otimes H^*$
to $H^*$ in $(^B_B{\cal YD}, C)$.

(iii) Let  $f$ be a $k$-linear map  from $U$ to $V$ and $g$
$k$-linear from $V$ to $W$ with $U, V, W$ in ${}^B_B{\cal YD }.$ If
$f$ and $gf$ are two morphisms in ${} ^B_B{\cal YD}$ with $Im (f ) =
V$, then $g$ is a morphism  $ {}^B_B{\cal YD }$.

(iv) Let  $M$ be an $H^*$-module in  ${} ^B_B{\cal YD}$, then $M$
has the maximal $H^*$-submodule $M^{rat}$ in ${} ^B_B{\cal YD}$, .

 \end {Lemma}

{\bf Proof.}
 (i) It is clear that $M^*$ is a $B$-module and $B$-comodule.
 For any $b \in B, h^* \in M^*, h\in M,$
on the one hand, $ \sum (b\cdot h^*)_{(-1)}<(b \cdot h^*)_{(0)}, h>
=\sum S^{-1}( h _{(-1)}) <h^*, S(b)\cdot h_{(0)}>.$ On the other
hand,
\begin {eqnarray*}
\sum  b_1 &h^*_{(-1)}S(b_3)&  <b_2\cdot h^*_{(0)},h> \\
&=& \sum   b _1 h^* _{(-1)} S(b_3) < h^*_{(0)}, S(b_2)\cdot h>\\
&=&\sum   b_1 S^{-1}((S(b_2)\cdot h)_{(-1)}) S(b_3)<h^*, (S(b_2)\cdot h)_{(0)}>\\
&=& \sum  S^{-1} (h_{(-1)}) b_2 S(b_3) <h^*, S(b_1)\cdot h_{(0)}>  \\
&=&\sum   S^{-1}(h _{(-1)}) <h^*, S(b)\cdot h_{(0)}>.\\
\end {eqnarray*}
Thus $M^*$ is a Yetter-Drinfeld $B$-module.

Obviously, $<,>$ is a $B$-module homomorphism. In order to show that
$<,>$ is a $B$-comodule homomorphism, it is enough to prove  that
$\sum h_{(-1) }^* h _{(-1)} <h^*_{(0)}, h_{(0)}> = 1_B <h^*, h>$ for
any $h^* \in M^*, h \in M.$ Indeed, the left side $= \sum
S^{-1}(h_{(-1) 2}) h _{(-1)1} <h^*, h_{(0)}> = 1 _B <h^*, h>.$ This
complete the proof.

(ii)  For any $b \in B, h, x \in H, h^* \in H^*, $ we see that
\begin {eqnarray*}
< b \cdot (h \rightharpoonup h ^*), x>
&=&  < (h\rightharpoonup h^* ),S(b)\cdot x>\\
&=& <h^*, (S(b)\cdot x)h> \hbox {\ \ \ \ \ \ \ \ \ \ and }\\
\sum _b<(b_1 \cdot h) \rightharpoonup (b_2 \cdot h^*), x>
&=& <h^*, S(b_2) \cdot (x (b_1 \cdot h)) >\\
&=& < h^*, (S(b_2)_1 \cdot x)(S(b_2)_2\cdot ( b_1 \cdot h) ) > \\
&=&  < h^*, (S(b_3) \cdot x) ((S(b_2)b_1)\cdot h) >\\
&=& <h^*, (S(b)\cdot x)h>.
\end {eqnarray*}
This show that $\rightharpoonup $ is a $B$-module homomorphism.
Similarly, we can show that it is a $B$-comodule homomorphism.

(iii) For any $b \in B, u \in U$, since $g (b \cdot f(u)) = gf (b
\cdot u) =b \cdot ( gf(u) )$, we have that $g$ is a $B$-module
homomorphism. Similarly, $g$ is a $B$-comodule homomorphism.

(iv) It can be shown by usual proof (see \cite [Theorem 2.2.6 and
Corollary 2.1.19] {DNR01}  ) that every $H^*$-submodule and quotient
$H^*$-module of rational $H^*$-module are rational. The direct sum
of rational $H^*$-modules is a rational. Consequently, The maximal
rational $H^*$-module $M^{rat}$ is the sum of all rational
$H^*$-modules of $M$.  $\Box$

Every $B$-module category $(_B{\cal M}, C^R)$ determined by
quasitriangulr Hopf algebra $(B,R)$
 is a full subcategory of  Yetter-Drinfeld
module category $(^B_B{\cal YD}, C)$.   Indeed, for any $B$-module
$(V, \alpha   )$, define $\phi (v) = \sum R_i ^{(2)} \otimes
R^{(1)}_i \cdot v$  for any $v \in V$, where $R = \sum _i R_i ^{(1)}
\otimes R_i ^{(2)}$. It is easy to check that $(V, \alpha , \phi )$
is a Yetter-Drinfeld $B$-module. Similarly, every $B$-comodule
category $(^B{\cal M}, C^r)$ determined by coquasitriangulr Hopf
algebra $(B,r)$
 is a full subcategory of  Yetter-Drinfeld
module category $(^B_B{\cal YD}, C)$.

\begin {Example} \label {3.2} (Existence of integrals )
Let  $H$ be a braided Hopf algebra in  $ (^B_B{\cal YD}, C)$ and the
antipode $S$ of $B$  satisfy $ S= S^{-1}$. Then $$ H^{ * rat } \cong
\int _{H^*}^l \otimes H \ \ \ \ \ ( \hbox {as $H$-Hopf modules in }
(^B_B{\cal YD}, C). )$$

\end {Example}

\begin {Example} \label {3.3} (Existence of integrals ) Let $H$ be a braided Hopf algebra in $({}^B_B{\cal YD }, C).$
 If $\lambda $ is a non-zero integral of $H\# B$ with $\lambda (a \otimes b) \not=0 $ for some $a \in H, b\in B$,
 then $\lambda ( id \otimes b)$ is a non-zero integral of
$H$, where $\lambda (id \otimes b)$ denote the $k$-linear map from
$H$ to $k$ by sending $h$ to $\lambda (h \otimes b)$ for any $h\in
H.$


\end {Example}
{\bf Proof.} It follows from \cite [Theorem 9.4.12] {Ma95b} and
\cite [P11] {AS02} that the bosonization $H \# B$ of braided Hopf
$H$ is a Hopf algebra.  For any $f \in H^*$ and any $x\in H$, we see

\begin {eqnarray*}
(f * \lambda (id \otimes b ) )(x) &=&  \sum f (x_1)\lambda  (x_2 \otimes b)\\
&=& ((f \otimes \epsilon _B) *\lambda ) (x \otimes b) \\
&=& f(1) \lambda ( x \otimes b).
\end {eqnarray*}
Thus
 $\lambda ( id \otimes b)$ is a non-zero  integral of $H$. $\Box$

{\bf Remark:} In Example \ref {3.3} it is possible that $B$ is
infinite-dimensional.

\begin {Example}\label {3.5}(Existence and uniqueness of integrals)( see \cite [Example 9.4.9]{Ma95b})
 If $H$ is an ordinary  coquasitriangular
Hopf algebra with a non-zero integral $\lambda$, then the braided
group analogue  $\underline H$ of H  has a non-zero integral
$\lambda $ in braided tensor category $({}^H {\cal M},C^r)$.
Conversely, if $\underline H$ has a non-zero integral, then so does
$H$. Indeed, since the comultiplication operations  of $H$ and
$\underline H$ are the same, we have that the multiplications of
$H^*$ and $\underline H^*$ are the same, so the integrals  of $H$
and $\underline H$ are the same. $\Box$   \end {Example}

\begin {Example} \label {3.4}  (The uniqueness of integrals ) Let $(H, \alpha , \phi )$ be a braided Hopf algebra in $({}_B ^B{\cal YD}, C)$
and $B$ is a finite-dimensional Hopf algebra. If $\phi $ is trivial,
then $dim \ \int _{H^*}^l =0$ or $1$.
\end {Example}
{\bf Proof.} Assume that $H$ has two linearly independent  non-zero
integrals $u^*$ and $w^*$. Let $v^*$ is a non-zero integral of $B$.
By \cite [Lemma 1.3.2] {DNR01}, $(H \otimes B)^* = H^* \otimes B^*$
as vector spaces. Since the $B$-comodule operation of $H$ is
trivial, we have that  $u^* \otimes v^*$ and $w^* \otimes v^*$ are
two linear independent  integrals of $H\# B$. This contradicts to
the fact $dim \int _ {(H\#B)^*}^l = 0$ or 1 (see, \cite [Theorem
5.4.2]{DNR01}). $\Box$

\section {Maschke's theorem for braided Hopf algebras}

In this section we give the relation between the integrals and
semisimplicity of braided Hopf algebras. Although authors in \cite
{GG03} gave the Maschke's theorem for rigid braided Hopf algebras ,
it is not known if every semisimple braided Hopf algebra is rigid or
finite. Thus our research of  the Maschke's theorem for infinite
braided Hopf algebras is useful.

Throughout this section we assume that there exists a forgetful
functor $F: \ \ {\cal C} \rightarrow {}_k{\cal M }$,
 such that $F(U \otimes V) = F(U)\otimes F(V)$
and $F(I)= k$, where $k$ is a field.

\begin {Theorem} \label {4.1} (The Maschke's theorem)
If $H$ is a finite dimensional braided Hopf algebra living in a
braided tensor category ${\cal C}$, then
 $H$ is  semisimple as ordinary algebra over field $k$
iff $\epsilon (\int _H^l)\not=0;$
\end {Theorem}

{\bf Proof.}
 If $H$ is semisimple then there is a left ideal $I$
such that
$$H = I\oplus ker \epsilon. $$
For any $y\in I, h\in H$, we see that
\begin {eqnarray*}
 hy &=& ((h- \epsilon (h) 1_H ) + \epsilon (h)1_H )y
\\
&=& (h- \epsilon (h) 1_H )y + \epsilon (h) y  \\
  &=& \epsilon (h)y    \hbox { \ \ \ \ \ since }
(h-\epsilon (h) 1_H) y \in
  (ker \epsilon )I =0 \ .
\end {eqnarray*}
Thus $y \in \int _H^l$, and so  $I \subseteq \int _H^l,$ which
implies $\epsilon (\int _H^l)\not=0.$

Conversely, if $\epsilon (\int _H^l) \not=0,$

let $z \in \int _H^l$ with $\epsilon (z) =1.$

Say $M$ is a left $H$-module and $N$ is an $H$-submodule of $M$.
 Assume that $\xi$  is a $k$-linear projection from
$M$ to $N$.
 We define $$\mu (m) = \sum z_1 \cdot \xi (S(z_2)\cdot
m)$$
 for every  $m \in M.$ It is sufficient to show that
$\mu$ is an $H$-module projection
 from $M$ to $N$. Obviously, $\mu$  is a $k$-linear
projection. Now we only need to show that it is an $H$-module map.
We see that $\alpha(H\otimes \mu)=\alpha(H\otimes \alpha)(H\otimes
id\otimes \xi)(H\otimes id\otimes \alpha)(H\otimes id\otimes
id\otimes m\otimes m)(H\otimes id\otimes S\otimes S\otimes H\otimes
M)(H\otimes \Delta(z)\otimes \Delta)(\Delta\otimes M)=
\alpha(H\otimes \xi)(H\otimes \alpha)(H\otimes m\otimes M)(H\otimes
S\otimes H\otimes M)(H\otimes m\otimes H\otimes M)(m\otimes C\otimes
H\otimes M)(H\otimes \Delta(z)\otimes \Delta\otimes M)(\Delta\otimes
M)=\alpha(H\otimes \xi)(H\otimes \alpha)(H\otimes m\otimes
M)(H\otimes S\otimes H\otimes M)(m\otimes m\otimes H\otimes
M)(H\otimes C\otimes id\otimes H\otimes M)(\Delta\otimes
\Delta(z)\otimes H\otimes M)(\Delta\otimes M)=\alpha(H\otimes
\xi)(H\otimes \alpha)(H\otimes m\otimes M)(H\otimes S\otimes
H\otimes M)(\Delta\otimes H\otimes M)(m\otimes H\otimes M)(H\otimes
C\otimes M)(\Delta\otimes z\otimes M)=\alpha(id\otimes
\xi)(id\otimes \alpha)(id\otimes m\otimes M)(id\otimes S\otimes
H\otimes M)(\Delta(z)\otimes H\otimes M)=\mu\circ\alpha$. Thus $\mu
$ is an $H$-module morphism. $\Box$

{\bf Remark:}   Theorem \ref {4.1} needs not $C_{H,H} =
C_{H,H}^{-1}$.

It is well-known that an ordinary algebra $H$ over a field $k$ is
called semisimple if every $H$-submodule $N$ of every $H$-module $M$
is a direct summand, i.e. if there is a $H$-submodule $L$  such that
$M = N \oplus L$ . Similarly we have the following definition.
Algebra $H$ in ${\cal C}$ is called semisimple with respect to
${\cal C}$, if every $H$-submodule $N$ in ${\cal C}$ of every
$H$-module $M$ in ${\cal C}$ is a direct summand( i.e. there is a
$H$-submodule $L$ in ${\cal C}$ such that $M = N \oplus L$ ).

\begin {Theorem} \label {4.2}
Let $H$ be a braided Hopf algebra in ${\cal C}$. If $H$ is
semisimple with respect to ${\cal C}$ and $ker \epsilon \in {\cal
C}$, then $\epsilon (\int _H^l)\not=0.$
\end {Theorem}
{\bf Proof.} It is similar to the proof of Theorem \ref {4.1}.
$\Box$

 \begin {Example} \label {4.3}

(see \cite [ P510 ] {Ma95b} )  Let $H= {\bf C} [x]$ denote
  the braided  line algebra. It is just the usual
algebra ${\bf C} [x]$
  of polynomials in $x$ over complex field   ${\bf
C}$, but we regard it
  as a q-statistical Hopf algebra with

 $$  \Delta (x) = x \otimes 1 + 1 \otimes x ,  \ \ \
 \epsilon (x ) =0 , \ \ \  S(x ) = -x , \ \ \ \ \mid x
^n \mid = n $$
 and $$C^r (x^n , x^m) =q^{nm} (x^m \otimes x^n).$$
In fact, $H$ is a braided Hopf algebra in $(^{\bf CZ}    {\cal M} ,
C^r)$ with coquasitriangular $r(m,n) = q^{mn}$. Here $\mid x^n \mid
$ denote the degree of $x^n$. If $y = \sum _0^n a_i x^i \in
\int_H^l,$ then $a_i=0 $  for $i =0, 1, 2, \cdots , n$ since $xy =
\epsilon (x)y =0.$ Thus  $\int _H^l =0$. It follows from  Theorem
\ref {4.1}  that
 $H$ is not semisimple.
$\Box$

 \end {Example}

\vskip 0.5cm {\bf Acknowledgement }: The work is supported by The
council of Hunan educations.


\begin{thebibliography}{150}
\bibitem {AS02} N. Andruskewisch and H.J.Schneider, Pointed Hopf algebras,
new directions in Hopf algebras, edited by S,Montgomery and
H.J.Schneider, Cambradge University Press, 2002.

\bibitem {BKL00} Y.Bespalov, T. Kerler and V.Lyubashenko.
Integral for braided Hopf algebras, J. Pure and Applied Algebras,
{\bf 148}(2000), 113--164.
\bibitem {DNR01} S.Dascalescu, C.Nastasecu, S. Raianu,
Hopf algebras: an introduction,  Marcel Dekker Inc. , 2001.

 \bibitem {GG03} J.A. Guccione and J.J.Guccione, Theory of braided Hopf crossed products,  J. Algebra, {\bf 261}(2003),
54--101.
 \bibitem {JS86}  A.Joyal and R.Street, Braided
monoidal categories , Math.Reports 86008; Macquaries University,
1986.

\bibitem {Ke99}  T.Kerler, Bridged links and tangle presentations of
cobordism categories, Adv. Math. {\bf 141} (1999), 207-- 281.

\bibitem {Ku96}  G.Guperberg, Non-involutary Hopf algebras and 3-manifold
invariants, Duke Math. J. {\bf 84} (1996), 83--129.

\bibitem {Lu93} G.Lusztig, Introduction to    quantum groups,
Progress on Math., Vol.110, Birkhauser, 1993.
\bibitem {LS69}  R.G.Larson, M.E.Sweedler, An associative orthogonal bilinear
form for Hopf algebras, Amer. J. Math. {\bf XCI} (1969)75-94.

\bibitem {Ma93} S. Majid,  Braided groups. J. Pure and Applied Algebra,
 {\bf 86}  (1993), 187--221.

\bibitem {Ma94} S. Majid,  Cross products by braided groups and
 bosonization, J. Algebra, {\bf 165}(1994), 165--190.

\bibitem {Ma95a} S.Majid, Algebras and Hopf algebras in braided categories,
Lecture notes in pure and applied mathematics advances in Hopf
algebras, Vol. 158, edited by J.Bergen and S.Montgomery,   1996.

\bibitem {Ma95b} S.Majid, Foundations of  quantum
group theory,  Cambradge University Press, 1995.

 \bibitem {Su71} J.B.Sulivan, The Uniqueness of integrals for Hopf algebras
 and some existence theorem of integrals
 for commutative Hopf algebras, J. Algebra, 19 (1971), 426-440.

\bibitem {Sw69a}   M.E.Sweedler, Hopf algebras, Benjamin, New York, 1969.

 \bibitem {Sw69b}   M. E. Sweedler, Integrals for Hopf algebras. Ann. of Math.
{\bf 89} (1969), 323--335.

 \bibitem {Ta99}   M. Takeuchi. Finite Hopf algebras in braided tensor
 categories, J. Pure and Applied Algebra, {\bf 138}(1999), 59-82.

 \bibitem {Tu94}   V.Turaev, Quantum invariants of knots and 3-manifolds,
 de Qruyter, Berlin, 1994.

\bibitem {ZC01} Shouchuan Zhang and Hui-Xiang Chen, The
double bicrossproducts in braided tensor categories,
    Communications in algebra,   {\bf 29}(2001), 31--66.
\end {thebibliography}
\end {document}